\documentclass[amssymb,twoside,10pt,amscd,reqno]{amsart}

\usepackage{latexsym,amsfonts}
\usepackage{amsmath,amsthm}          
\usepackage{euscript}                

\theoremstyle{plain} \numberwithin{equation}{section}
\newtheorem{thm}{Theorem}[section]

\newtheorem{prop}[thm]{Proposition}

\newtheorem{lemma}[thm]{Lemma}

\newtheorem{claim}[thm]{Claim}
\newtheorem{conj}[thm]{Conjecture}

\theoremstyle{definition}

\newtheorem{defn}[thm]{Definition}

\newtheorem{rmk}[thm]{Remark}

\def\CC{\mathbb{C}}
\def\QQ{\mathbb{Q}}
\def\PP{\mathbb{P}}
\def\ZZ{\mathbb{Z}}

\def\ka{\kappa}

\def\Th{\Theta}

\def\Lb{\mathcal{L}}
\def\E{\mathcal{E}}
\def\O{\mathcal{O}}
\def\Ext{\text{Ext}}
\def\h{\text{h}}
\def\H{\text{H}}
\def\Hom{\text{Hom}}
\def\Pic{\text{Pic}}
\def\Sec{\text{Sec}}
\def\Sym{\text{Sym}}

\def\dra{\dashrightarrow}
\def\ra{\rightarrow}
\def\setm{\setminus}

\begin{document}

\title{Examples of Fano varieties of index one that are not birationally rigid}
\author[Castravet]{Ana-Maria Castravet}
\address{Mathematics, Univ. of Texas at Austin\\ Austin, TX 78712}
\email{noni@math.utexas.edu}
\date{\today}

\maketitle

\section{Pukhlikov's Conjecture}

The following is a conjecture of Pukhlikov:

\begin{conj}\cite[Conjecture~5.1]{P4}\label{main}
Let $V$ be a smooth Fano variety of dimension $\dim V\geq4$ with
Picard group $\Pic(V)$ generated by the canonical class ${K_V}$.
Then $V$ is birationally rigid. If $\dim V\geq5$, then $V$ is
superrigid.
\end{conj}

Conjecture \ref{main} was proved for a large class of Fano complete
intersections \cite{P1}, \cite{P3}, \cite{P4}, \cite{dF},
\cite{dFEM}. More examples are given by complete intersections in
weighted projective spaces \cite{P2}. Conjecture \ref{main} is a
generalization of the famous theorem of Iskovskikh and Manin
\cite{IM} that states that a smooth quartic threefold is
birationally rigid, and therefore not rational.

\

In this note we give counterexamples (in arbitrarily large
dimension) to Conjecture \ref{main} using moduli spaces of bundles
on curves. I am grateful to Professor J\'anos Koll\'ar for telling
me about this question and suggesting to look at sections of the
theta divisor on moduli spaces of bundles on curves. I thank Jenia
Tevelev and Sean Keel for reading this note and for helpful
suggestions.  

\

We recall the basic definitions about birational rigidity from
\cite{P4}. Let $V$ be a uniruled $\QQ$-Gorenstein variety with
terminal singularities. For an effective divisor $D\neq0$ on $V$,
one defines the \emph{canonical threshold of canonical adjunction}
$c(D)$ of the divisor $D$ as follows:

$$c(D)=\hbox{sup }\{\quad b/a \quad|\quad
b,a\in\ZZ_+\setminus\{0\},\quad |aD+bK_X|\neq\emptyset\}.$$

(If there are no $a, b\in\ZZ_+\setminus\{0\}$ such that
$|aD+bK_X|\neq\emptyset$, then we set $c(D)=0$.)

\

Note, for an effective divisor $D\neq0$ on a Fano variety $V$ with
$\Pic(V)=\ZZ{K_V}$, the canonical threshold $c(D)$ is the number
$-m$, if $D={m}K_V$ in $\Pic(V)$. Clearly, in this case $c(D)>0$.

\

In what follows all varieties are assumed to be $\QQ$-factorial with
terminal singularities. For simplicity, we work over an
algebraically closed field of characteristic zero.

\begin{defn}\cite[Def. 5.1]{P4}\label{def1}
A variety $V$ is called \emph{birationally rigid}, if for any $V'$,
any birational map $\phi:V\dra V'$ and any moving linear system
$\Sigma'$ on $V'$, there exists a birational self-map $\alpha:V\dra
V$ such that if $\Sigma$ is the birational transform of $\Sigma'$
via the composition $\phi\circ\alpha$ (i.e., the linear system
induced on $V$ when composing with $\phi\circ\alpha$), one has
$c(\Sigma)\leq c(\Sigma')$. The variety $V$ is called
\emph{birationally superrigid}, if one may always take $\alpha=id$.
\end{defn}

The following Proposition is an immediate consequence of the above
definitions:

\begin{prop}\cite[Prop.~5.1]{P4}\label{fibrations}
Let $V$ be a smooth Fano variety with $\Pic(V)=\ZZ{K_V}$. If $V$ is
birationally rigid, then it is impossible to have a birational map
$V\dra V'$, with $V'\ra S'$ a morphism with uniruled general fiber
and $S'$ a projective variety of dimension $\dim S'\geq1$.
\end{prop}

\begin{proof} Assume there is a birational map $V\dra V'$, with
$\pi:V'\ra S'$ a morphism with uniruled general fiber and $\dim
S'\geq1$. Let $D'$ be an effective Cartier divisor on $S'$ and let
$\Sigma'=\pi^*D'$. By Definition \ref{def1}, there is a birational
map $\alpha:V\dra V$ such that if $\Sigma$ is the birational
transform of $\Sigma'$ via the composition $\phi\circ\alpha$, then
$c(\Sigma)\leq c(\Sigma')$. Note that since $V$ is a Fano variety
with $\Pic(V)=\ZZ{K_V}$, one has $c(\Sigma)>0$.

\

We claim that $c(\Sigma')=0$: if there are $a, b>0$ such that
$a\pi^*D'+bK_{V'}$ is effective on $V'$, then for a general fiber
$F$ of $\pi$ (choose $F$ uniruled, in the smooth locus of $\pi$ and
such that it is not contained in the divisor $a\pi^*D'+bK_{V'}$),
one has that the divisor
$$(a\pi^*D'+bK_{V'})_{|F}=(bK_{V'})_{|F}=bK_F$$
is effective on $F$. This is a contradiction, since $F$ uniruled
implies that $\H^0(bK_F)=0$ for all $b>0$ \cite{K}.
\end{proof}

Note that since the condition of being uniruled is a closed
condition \cite{K}, one may drop the word ``general" from the
statement of Proposition \ref{fibrations}. Moreover, since the
condition of being uniruled is a birational property, one may as
well replace the condition of having a morphism $V'\ra S'$ with
having  a birational map $V'\dra S'$ with the same properties.

\

To motivate geometrically Definition \ref{def1}, we recall the other
definition of birational rigidity from \cite{Co}. Recall that a
morphism $f:V\ra S$ is called a \emph{Mori fiber space} if it has
relative Picard number $1$, $-K_V$ is relatively ample for $f$ and
$\dim S<\dim V$. Note, Fano varieties with Picard number $1$ are
trivially Mori fiber spaces. A birational map $\phi:V\dra V'$ to
another Mori fiber space $f':V'\ra S'$ is \emph{square} if there is
a birational map $h:S\dra S'$ such that $h\circ f=f'\circ\phi$ and
the map $\phi$ induces an isomorphism on the general fibers.

\begin{defn}\cite[Def. 1.3]{Co}\label{def2}
A Mori fiber space $f:V\ra S$ is called \emph{birationally rigid}
if for any birational map $\phi:V\dra V'$ to another Mori fiber
space $f':V'\ra S'$, there is a birational self-map $\alpha:V\dra
V$ such that $\phi\circ\alpha$ is square. If for any $\phi$ as
above it follows that $\phi$ is square, then $V$ is called
\emph{birationally superrigid}.
\end{defn}

A Mori fiber space that satisfies Definition \ref{def1} also
satisfies Definition \ref{def2} (by the Sarkisov program, based on
the Noether-Fano-Iskovskikh inequalities; this is known in
dimension $3$ and conjectured in higher dimensions) . By the Mori
program, any uniruled variety is birational to a Mori fiber space;
hence, if in addition one assumes the Mori program, the two
definitions are equivalent.

\

It follows from Definition \ref{def2} that if $V$ is a smooth Fano
variety with Picard number $1$ and birationally rigid, then there is
no birational map $\phi:V\dra V'$ with $f':V'\ra S'$ another Mori
fiber space with $\dim S'>0$. From the Mori program for the relative
case, one deduces Prop. \ref{fibrations}. In particular, note that
if $V$ is a smooth Fano variety with Picard number $1$ and
birationally rigid, then for any $V'$ another Fano variety of Picard
number $1$, if $V'$ is birational to $V$, then  $V\cong V'$.

\section{Counterexamples using moduli spaces of bundles on curves}

Let $C$ be a smooth projective curve over $\CC$ of genus $g\geq3$.
Fix $\xi$ to be a degree $1$ line bundle on $C$. Let $M$ be the
moduli space of stable, rank $2$ vector bundles on $C$ with
determinant $\xi$. The moduli space $M$ is a smooth, projective,
variety of dimension $3g-3$. The Picard group of $M$ is $\ZZ$
\cite{DN}. Let $\Th$ be the ample generator. In fact, $\Th$ is
very ample \cite{BV}. Then $K_M\cong-2\Th$ \cite{R}.

\smallskip

Let $N$ be a nonsingular element of the linear system $|\Th|$. Let
$\Th'\in\Pic(N)$ be the restriction of $\Th$ to $N$. The canonical
bundle of $N$ is $-\Th'$. Since $g\geq3$, by Lefschetz's theorem,
$\Th'$ generates $\Pic(N)$. Therefore, the variety $N$ satisfies
the conditions in Conjecture \ref{main}. We prove the following:

\begin{prop}\label{example}
If $N$ is a general element of the linear system $|\Th|$,
then $N$ is a smooth Fano variety with $\Pic(N)\cong\ZZ K_N$ that
is not birationally rigid.
\end{prop}

\begin{proof}
We start with a general construction. Let $e\geq0$ and let $\Lb$
be a line bundle of degree $-e$ on $C$. Denote by $V_{\Lb}$ the
space of extensions
$\Ext^1(\Lb^{-1}\otimes\xi,\Lb)\cong\H^1(C,\Lb^2\otimes\xi^{-1})$.
Then $V_{\Lb}$ parametrizes extensions of the form
\begin{equation}\tag{$*$}\label{*}
0\ra\Lb\ra\E\ra\Lb^{-1}\otimes\xi\ra0
\end{equation}
The dimension of $V_{\Lb}$ is $2e+g$. Clearly, any two non-zero
elements in $V_{\Lb}$ which differ by a scalar define isomorphic
vector bundles $\E$. Therefore, the isomorphism classes of
non-trivial extensions as above are parametrized by the projective
space $\PP(V_{\Lb})$. The locus $Z_{\Lb}\subset\PP(V_{\Lb})$ of
extensions (\ref{*}) with $\E$ unstable is an irreducible
subvariety of codimension at least $g$ \cite[Lemma~2.1]{Ca}. There
is a well-defined morphism
\begin{equation}\label{ka}
\ka_{\Lb}:\PP(V_{\Lb})\setm Z_{\Lb}\ra M
\end{equation}
that associates to an extension (\ref{*}) the isomorphism class of
$\E$. By \cite[2.3(1), or Lemma~A.1]{Ca} $$\ka_{\Lb}^*\Th\cong
\O(2e+1).$$

Consider the case when $e=g-1$. Then $\PP(V_{\Lb})\cong\PP^{3g-3}$
and $\ka_{\Lb}^*\Th\cong\O(2g-1)$. The following Claim is a
standard fact. For convenience, we include a short proof.

\begin{claim}\label{biratl} If $\Lb$ is a line bundle of degree $1-g$,
 the morphism $\ka_{\Lb}$ in (\ref{ka}) is birational.
\end{claim}

\emph{Proof of Claim \ref{biratl}. } Note that for any $\E\in M$,
by Riemann-Roch one has
$$\chi(\E^*\otimes\Lb^{-1}\otimes\xi)=1.$$

Hence, for any $\E\in M$ one has:
\begin{equation*}
\Hom(\E,\Lb^{-1}\otimes\xi)\cong\H^0(\E^*\otimes\Lb^{-1}\otimes\xi)\neq0
\end{equation*}


If for general $\E\in M$, any non-zero morphism
$\phi:\E\ra\Lb^{-1}\otimes\xi$ is  surjective, then we are done,
as $\phi$ determines uniquely (up to scaling)
the extension (\ref{*}), and by dimension considerations, one must have
that $$\h^0(\E^*\otimes\Lb^{-1}\otimes\xi)=1,$$
(i.e., the fiber of $\ka$ at a general point contains a unique closed point).

We prove that for general $\E$, a non-zero morphism
$\phi:\E\ra\Lb^{-1}\otimes\xi$ must be surjective. Note that if
$\phi$ is not surjective, then its image is
$\Lb^{-1}\otimes\xi(-D)$, for some effective divisor $D$ of degree
$d>0$ and there is an exact sequence:
\begin{equation}\label{extension}
0\ra\Lb(D)\ra\E\ra\Lb^{-1}\otimes\xi(-D)\ra0
\end{equation}

It follows from the stability of $\E$ that $d<g$. For each $0<d<g$
construct the total space $\PP_d$ of extensions (\ref{extension})
by letting $D$ vary in $\Sym^d(C)$: the space $\PP_d$ is a
projective bundle over $\Sym^d(C)$ with fiber at $D$ isomorphic to
$\PP(V_{\Lb(D)})$. The dimension is:
$$\dim\PP_d=d+\dim V_{\Lb(D)}-1=d+2(g-1-d)+g-1=3g-3-d.$$

The vector bundle which is the middle term of the universal
extension over $\PP_d$ induces a rational map $\rho:\PP_d\dra M$.
Since $d>0$, the map $\rho$ is not dominant. Therefore, a general
$\E\in M$ will not sit in an exact sequence (\ref{extension}).

\

From the previous discussion and since $\Th$ is very ample, one
has the following:
\begin{claim}\label{hypers} If $N$ is a general element of the linear
system $|\Th|$, then $N$ is birational to an irreducible
hypersurface $X_{2g-1}$ in $\PP^{3g-3}$ of degree $2g-1$.
\end{claim}

Proposition \ref{example} follows now from Claim \ref{hypers} and
Lemma \ref{hypers of small degree}.
\end{proof}

\begin{lemma}\label{hypers of small degree}
For any irreducible hypersurface $X_d\subset\PP^n$ of degree $d<n$
(possibly singular) there is a rational map $\rho:X\dra S$, $\dim
S>0$, with uniruled fibers.
\end{lemma}

\begin{proof} Let $S$ be a general pencil of hyperplanes $\{H_s\}_{s\in
S}$ in $\PP^n$. Then $X_s=X\cap H_s$ are hypersurfaces of degree
$d$ in $H_s\cong\PP^{n-1}$. A smooth hypersurface of degree $d<n$
in $\PP^{n-1}$ is Fano; hence, it is rationally connected
\cite{KMM}. It follows by deformation theory that any irreducible
hypersurface of degree $d<n$ in $\PP^{n-1}$ is uniruled. Hence,
for all $s\in S$ such that $X_s$ is irreducible, $X_s$ is
uniruled. Therefore, the induced rational map $\rho:X\dra S$ has
uniruled fibers.
\end{proof}


\begin{rmk}
The cohomology group $\H^4(M;\QQ)$ has two independent generators
\cite{N}. Hence, by Lefschetz's theorem, if $g\geq4$, the rank of
$\H^4(N;\QQ)$ is also $2$. Since again by Lefschetz's theorem, the
cohomology group $\H^4$ of any smooth complete intersection of
dimension $\geq5$ is of rank $1$, it follows that $N$ is not a
complete intersection.
\end{rmk}

\section{Description of the hypersurface $X_{2g-1}\subset\PP^{3g-3}$}

Our construction of the map $\ka_{\Lb}$ is a variant of the
construction of Bertram \cite{B}.  The above results about the map
$\ka_{\Lb}$ ($\deg\Lb=1-g$) also follow from \cite{B}. We chose to
include the above considerations (which are enough for the purpose
of this note) because of their simplicity and to avoid referring
to the technical results in \cite{B}. However, Bertram's powerful
construction gives a precise description of the hypersurface
$X_{2g-1}\subset\PP^{3g-3}$. We describe this below.

\

Let $\Lb$ be a line bundle of degree $1-g$ on $C$. Then $C$ has a
natural embedding $C\subset\PP(V_{\Lb})\cong\PP^{3g-3}$ given by
$\Lb^{-2}\otimes\xi\otimes K_C$, since by Serre duality one has:
$$V_{\Lb}\cong\H^1(C,\Lb^2\otimes\xi^{-1})\cong
\H^0(C,\Lb^{-2}\otimes\xi\otimes K_C)^*.$$

Let $\Sec^k(C)$ be the $(k+1)$-secant variety of $C$ (i.e., the
closure in $\PP^{3g-3}$ of the union of all the $k$-planes spanned
by $k+1$ distinct points on $C$).

\begin{thm}\cite[Thm. 1]{B}\label{bertram1}
There is a sequence of blow-ups
$\pi:\tilde{\PP}_{\Lb}\ra\PP(V_{\Lb})$ with smooth centers
(starting with the blow-up of $\PP(V_{\Lb})$ along $C$) that
resolves the rational map $\ka_{\Lb}:\PP(V_{\Lb})\dra M$ into a
morphism $\tilde{\PP}_{\Lb}\ra M$. There are $g$ exceptional
divisors $E_0,E_1,\ldots,E_{g-1}$ and  $E_k$ dominates the secant
variety $\Sec^k(C)$ for every $k$.
\end{thm}

\begin{thm}\cite[Thm. 2, Prop. 4.7]{B}\label{bertram2}
There is a natural identification
\begin{equation}\label{lin system}
\H^0(M,\Th)\cong\H^0(\tilde{\PP}_{\Lb},(2g-1)H-(2g-3)E_0-(2g-5)E_1-\ldots-E_{g-2}).
\end{equation}
\end{thm}

It follows from Theorem \ref{bertram2} that the proper transform
$\tilde{X}$ in $\tilde{\PP}_{\Lb}$ of the hypersurface
$X_{2g-1}\subset\PP^{3g-3}$ of Proposition \ref{hypers} is a
general member of the linear system in (\ref{lin system}). Hence,
by Bertini, $\tilde{X}$ is smooth.

\begin{prop}
The singular locus of the hypersurface $X_{2g-1}\subset\PP^{3g-3}$
has codimension $\geq g-2$.  Hence, if $g\geq4$ then $X$ is normal
and the canonical bundle $K_X$ is Cartier. Moreover, in this case
$X$ has terminal singularities.
\end{prop}

\begin{proof}
Since the proper transform $\tilde{X}$ of $X$ is smooth, it
follows that $X$ is smooth outside $\Sec^{g-1}(C)\cap X$. Since a
general $X$ does not contain $\Sec^{g-1}(C)$, it follows that the
singular locus of $X$ has dimension at most
$\dim\Sec^{g-1}(C)-1=2g-2$. It is well-known that hypersurfaces
$X$ whose singular locus has codimension at least $2$ are normal
and the canonical class $K_X$ is Cartier. Hence, if $g\geq4$ then
$X$ is normal and has the canonical class:
$$K_X=\O_X(1-g).$$

Consider the resolution $\pi:\tilde{X}\ra X$. The canonical class
of $\tilde{\PP}_{\Lb}$ is given by:

$$K_{\tilde{\PP}_{\Lb}}=-(3g-2)H+(3g-5)E_0+(3g-7)E_1+\ldots+(g-1)E_{g-2}$$

\

The canonical class of $\tilde{X}$ is given by:

$$K_{\tilde{X}}=(K_{\tilde{\PP}_{\Lb}}+\tilde{X})_{|\tilde{X}}=
-(g-1)H+(g-2)E_0+(g-2)E_1+\ldots+(g-2)E_{g-2}.$$

Hence, $X$ has terminal singularities.

\end{proof}

\end{document}